\input amstex

\input amstex
\input amsppt.sty
\magnification=\magstep1
\hsize=36truecc
 \vsize=23.5truecm
\baselineskip=14truept
 \NoBlackBoxes
\def\q{\quad}
\def\qq{\qquad}
\def\mod#1{\ (\text{\rm mod}\ #1)}

\def\t{\text}
\def\qtq#1{\q\t{#1}\q}
\def\mod#1{\ (\text{\rm mod}\ #1)}
\def\qtq#1{\q\t{#1}\q}
\def\f{\frac}
\def\e{\equiv}
\def\b{\binom}

\def\sls#1#2{(\f{#1}{#2})}
 \def\ls#1#2{\big(\f{#1}{#2}\big)}
\def\Ls#1#2{\Big(\f{#1}{#2}\Big)}
%\par\q The paper is submitted for possible publication in
%Duke Mathematical Journal.
%\par\q  April 28, 2011
\par\q
\let \pro=\proclaim
\let \endpro=\endproclaim
\topmatter
\title Congruences involving $\b{2k}k^2\b{3k}km^{-k}$
\endtitle
\author ZHI-Hong Sun\endauthor
\affil School of the Mathematical Sciences, Huaiyin Normal
University,
\\ Huaian, Jiangsu 223001, PR China
\\ Email: zhihongsun$\@$yahoo.com
\\ Homepage: http://www.hytc.edu.cn/xsjl/szh
\endaffil

 \nologo \NoRunningHeads

\abstract{Let $p>3$ be a prime, and let $m$ be an integer with
$p\nmid m$. In the paper, based on the work of Brillhart and Morton,
by using the work of Ishii and Deuring's theorem for elliptic curves
with complex multiplication we solve some conjectures of Zhi-Wei Sun
concerning $\sum_{k=0}^{p-1}\binom{2k}k^2\binom{3k}km^{-k}\mod
{p^2}$.
  \par\q
 }
 \endabstract
  \footnote"" {The author is
supported by the Natural Sciences Foundation of China (grant No.
10971078).}
\endtopmatter
\document
\subheading{1. Introduction}
\par For positive integers $a,b$ and $n$, if $n=ax^2+by^2$
for some integers $x$ and $y$, we briefly say that $n=ax^2+by^2$.
Let $p>3$ be a prime. In 2003, Rodriguez-Villegas[RV] posed some
conjectures on supercongruences modulo $p^2$. One of his conjectures
is equivalent to
$$\sum_{k=0}^{p-1}\f{\b{2k}k^2\b{3k}k}{108^k}\e
\cases 4x^2-2p\mod{p^2}&\t{if $p=x^2+3y^2\e 1\mod 3$,}
\\0\mod{p^2}&\t{if $p\e 2\mod 3$.}
\endcases$$
  This conjecture has been solved  by
Mortenson[Mo] and Zhi-Wei Sun[Su2].
\par  Let $\Bbb Z$ be the set of
integers, and for a prime $p$ let $\Bbb Z_p$ be the set of rational
numbers whose denominator is coprime to $p$. Recently the author's
brother Zhi-Wei Sun[Su1] posed many conjectures for
$\sum_{k=0}^{p-1}\b{2k}k^2\b{3k}km^{-k}\mod{p^2}$,
 where $p>3$ is a prime and
$m\in\Bbb Z$ with $p\nmid m$. For example, he conjectured that (see
[Su1, Conjecture A13])
$$\aligned\sum_{k=0}^{p-1}\f{\b{2k}k^2\b{3k}k}{(-27)^k}\e
\cases 0\mod {p^2}&\t{if $p\e 7,11,13,14\mod {15}$,}
\\4x^2-2p\mod {p^2}&\t{if $p=x^2+15y^2\e 1,4\mod {15}$,}
\\20x^2-2p\mod {p^2}&\t{if $p=5x^2+3y^2\e 2,8\mod {15}$.}
\endcases\endaligned\tag 1.1$$

\par Let $\{P_n(x)\}$ be the Legendre polynomials given by
 (see [MOS, pp.\;228-232], [G,
(3.132)-(3.133)])
$$P_n(x)=\f
1{2^n}\sum_{k=0}^{[n/2]}\b nk(-1)^k\b{2n-2k}nx^{n-2k} =\f 1{2^n\cdot
n!}\cdot\f{d^n}{dx^n}(x^2-1)^n,\tag 1.2$$ where $[a]$ is the
greatest integer not exceeding $a$.
 From (1.2)
we see that
$$P_n(-x)=(-1)^nP_n(x).\tag 1.3$$
\par Let $\sls am$ be the Jacobi symbol. For a prime $p>3$, in [S2] the author showed
that
$$P_{[\f p3]}(t)
\e \sum_{k=0}^{[p/3]}\b{2k}k\b{3k}k\Ls{1-t}{54}^k\e
\sum_{k=0}^{p-1}\b{2k}k\b{3k}k\Ls{1-t}{54}^k\mod p.\tag 1.4$$ We
note that $p\mid\b{2k}k\b{3k}k$ for $\f p3<k<p$. In the paper, using
the work of Brillhart and Morton[BM] we prove that
$$P_{[\f
p3]}(t) \e -\Ls
p3\sum_{x=0}^{p-1}\Ls{x^3+3(4t-5)x+2(2t^2-14t+11)}p\mod p.\tag 1.5$$
Based on (1.5) and the work of Ishii[I], we determine
$$P_{[\f
p3]}(t)\mod p\qtq{for}t=\f 54,\f 5{\sqrt{-2}},\f{\sqrt{-11}}4,\f
1{\sqrt 2},\sqrt 5,\f{9}{20}\sqrt 5,\f{\sqrt{17}}4,\f
5{32}\sqrt{41},\f{53}{500}\sqrt{89}.$$ For instance, if $p\e 1,4\mod
5$ is a prime, then
$$P_{[\f p3]}(\sqrt 5)\e \cases
2x\sls x3\mod p&\t{if $p=x^2+15y^2\e 1,4\mod {15}$,} \\0\mod p&\t{if
$p\e 11,14\mod{15}$}.\endcases$$ \par Let $p>3$ be a prime,
$m\in\Bbb Z_p$, $m\not\e 0\mod p$ and $t=\sqrt{1-108/m}$. In the
paper we show that
$$\sum_{k=0}^{p-1}\f{\b{2k}k^2\b{3k}k}{m^k}\e
P_{[\f p3]}(t)^2\mod p\tag 1.6$$ and that
$$P_{[\f p3]}(t)\e 0\mod p\qtq{implies}
\sum_{k=0}^{p-1}\f{\b{2k}k^2\b{3k}k}{m^k}\e 0\mod{p^2}.\tag 1.7$$ On
the basis of (1.6) and (1.7), we prove some congruences for
$\sum_{k=0}^{p-1}\b{2k}k^2\b{3k}km^{-k}$ in the cases
$m=8,64,216,-27,-192,-8640,-12^3,-48^3,-300^3$. Thus we partially
solve some conjectures posed by Zhi-Wei Sun in [Su1]. As two
examples, for odd primes $p\not=11$ we have
$$\sum_{k=0}^{p-1}\f{\b{2k}k^2\b{3k}k}{64^k}\e\cases
 x^2\mod p&\t{if $\sls
p{11}=1$ and so $4p=x^2+11y^2$,}\\0\mod {p^2}&\t{if $\sls
p{11}=-1$,}\endcases$$ for odd primes $p$ with $\sls{17}p=1$ we have
$$\sum_{k=0}^{p-1}\f{\b{2k}k^2\b{3k}k}{(-12)^{3k}}
\e\cases x^2\mod p&\t{if $p\e 1\mod 3$ and so $4p=x^2+51y^2$,}
\\0\mod{p^2}&\t{if $p\e 2\mod 3$.}
\endcases$$

 \subheading{2. A general congruence modulo
$p^2$} \pro{Lemma 2.1} Let $m$ be a nonnegative integer. Then
$$\sum_{k=0}^m\b {2k}k^2\b{3k}k\b k{m-k}(-27)^{m-k}
=\sum_{k=0}^m\b {2k}k\b{3k}k\b{2(m-k)}{m-k}\b{3(m-k)}{m-k}.$$
\endpro We prove the lemma by using WZ method and Mathematica.
Clearly the result is true for $m=0,1$. Since both sides satisfy the
same recurrence relation
$$81(m+1)(3m+2)(3m+4)S(m)
    -3(2m+3)(9m^2+27m+22)S(m+1) + (m+2)^3 S(m+2) = 0,$$
    we see that the lemma is true.
The  proof certificate for the left hand side is
$$ - \frac{729 k^2(m+2)(m-2k)(m-2k+1)}{(m-k+1)(m-k+2)},
$$ and the proof certificate for the right hand side is $$
\frac{9
k^2(3m-3k+1)(3m-3k+2)(9m^2-9mk+30m-14k+24)}{(m-k+1)^2(m-k+2)^2}. $$
\pro{Theorem 2.1} Let $p$ be an odd prime and let  $x$ be a
variable. Then
$$\sum_{k=0}^{p-1}\b{2k}k^2\b{3k}k(x(1-27x))^k\e \Big(
\sum_{k=0}^{p-1}\b{2k}k\b{3k}kx^k\Big)^2\mod {p^2}.$$
\endpro
Proof. It is clear that
$$\align &\sum_{k=0}^{p-1}\b{2k}k^2\b{3k}k(x(1-27x))^k
\\&=\sum_{k=0}^{p-1}\b{2k}k^2\b{3k}kx^k\sum_{r=0}^k\b kr(-27x)^r
\\&=\sum_{m=0}^{2(p-1)}x^m\sum_{k=0}^{min\{m,p-1\}}\b{2k}k^2\b{3k}k\b
k{m-k}(-27)^{m-k}.\endalign$$
 Suppose $p\le m\le 2p-2$ and $0\le
k\le p-1$. If $k>\f p2$, then $p\mid \b{2k}k$ and so $p^2\mid
\b{2k}k^2$. If $k<\f p2$, then $m-k\ge p-k>k$ and so $\b k{m-k}=0$.
Thus, from the above and Lemma 2.1 we deduce that
$$\align &\sum_{k=0}^{p-1}\b{2k}k^2\b{3k}k(x(1-27x))^k
\\&\e \sum_{m=0}^{p-1}x^m\sum_{k=0}^m\b {2k}k^2\b{3k}k
\b k{m-k}(-27)^{m-k}
\\&=\sum_{m=0}^{p-1}x^m\sum_{k=0}^m\b {2k}k\b{3k}k
\b{2(m-k)}{m-k}\b{3(m-k)}{m-k}
\\&=\sum_{k=0}^{p-1}\b{2k}k\b{3k}kx^k\sum_{m=k}^{p-1}
\b{2(m-k)}{m-k}\b{3(m-k)}{m-k}x^{m-k}
\\&=\sum_{k=0}^{p-1}\b{2k}k\b{3k}kx^k\sum_{r=0}^{p-1-k}\b{2r}r\b{3r}rx^r
\\&=\sum_{k=0}^{p-1}\b{2k}k\b{3k}kx^k\Big(\sum_{r=0}^{p-1}\b{2r}r\b{3r}rx^r
-\sum_{r=p-k}^{p-1}\b{2r}r\b{3r}rx^r\Big)
\\&=\Big(\sum_{k=0}^{p-1}\b{2k}k\b{3k}kx^k\Big)^2
-\sum_{k=0}^{p-1}\b{2k}k\b{3k}kx^k
\sum_{r=p-k}^{p-1}\b{2r}r\b{3r}rx^r
 \mod{p^2}.\endalign$$
If $\f{2p}3\le k\le p-1$, then $\b{2k}k\b{3k}k=\f{(3k)!}{k!^3}\e
0\mod {p^2}$. If $0\le k \le \f p3$ and $p-k\le r\le p-1$, then
$\f{2p}3\le r\le p-1$ and so $\b{2r}r\b{3r}r=\f{(3r)!}{r!^3}\e 0\mod
{p^2}$. If $\f p3<k<\f{2p}3$ and $p-k\le r\le p-1$, then $r\ge
p-k>\f p3$, $\b{2k}k\b{3k}k=\f{(3k)!}{k!^3}\e 0\mod p$ and
$\b{2r}r\b{3r}r=\f{(3r)!}{r!^3}\e 0\mod p$. Hence, for $0\le k\le
p-1$ and $p-k\le r\le p-1$ we have $p^2\mid
\b{2k}k\b{3k}k\b{2r}r\b{3r}r$ and so
$$\sum_{k=0}^{p-1}\b{2k}k\b{3k}kx^k
\sum_{r=p-k}^{p-1}\b{2r}r\b{3r}rx^r
 \e 0\mod{p^2}.$$
Therefore the result follows.
\pro{Corollary 2.1} Let $p>3$ be a
prime and $m\in\Bbb Z_p$ with $m\not\e 0\mod p$. Then
$$\sum_{k=0}^{p-1}\f{\b{2k}k^2\b{3k}k}{m^k}
\e\Big(\sum_{k=0}^{p-1}\b{2k}k\b{3k}k\Big(\f{1-\sqrt{1-108/m}}{54}\Big)^k\Big)^2
\mod{p^2}.$$
\endpro
Proof. Taking $x=\f{1-\sqrt{1-108/m}}{54}$ in Theorem 2.1 we deduce
the result.

\subheading{3. Congruences for $P_{[p/3]}(t)\mod p$}
\par Let  $W_n(x)$ be the Deuring polynomial given by
$$W_n(x)=\sum_{k=0}^n\b nk^2x^k.\tag 3.1$$
It is known that ([G,(3.134)],[BM])
$$W_n(x)=(1-x)^nP_n\Ls{1+x}{1-x}.\tag 3.2$$
Let $p>3$ be a prime , $m,n\in\Bbb Z_p$ and $4m^3+27n^2\not\e 0\mod
p$.  From [Mor, Theorem 3.3] we have
$$\aligned&\sum_{x=0}^{p-1}\Ls{x^3+mx+n}p\\&\e -(-48m)^{\f{1-\sls p3}2}(864n)^{\f{1-\sls{-1}p}2}
(-16(4m^3+27n^2))^{[\f p{12}]}J_p\Ls{2^8\cdot
3^3m^3}{4m^3+27n^2}\mod p,\endaligned\tag 3.3$$ where $J_p(t)$ is a
certain Jacobi polynomial given by
$$J_p(t)=1728^{[\f p{12}]}P_{[\f p{12}]}^{(-\f 13\sls p3,-\f
12\sls{-1}p)}\Big(1-\f t{864}\Big)\tag 3.4$$ and
$$P_k^{(\alpha,\beta)}(x)=\f
1{2^k}\sum_{r=0}^k\b{k+\alpha}r\b{k+\beta}{k-r}(x-1)^{k-r}(x+1)^r.$$
 \pro{Theorem 3.1} Let $p>3$ be a
prime and $t\in\Bbb Z_p$. Then
$$P_{[\f p3]}(t)\e -\Ls p3
\sum_{x=0}^{p-1}\Ls{x^3+3(4t-5)x+2(2t^2-14t+11)}p\mod p.$$
\endpro
Proof. It is well known that $P_n(1)=1$. Since $P_{[\f p3]}(1)=1$
and
$$\sum_{x=0}^{p-1}\Ls{x^3-3x-2}p=\sum_{x=0}^{p-1}\Ls{(x+1)^2(x-2)}p
=\sum_{x=0}^{p-1}\Ls{x-2}p-\Ls{-1-2}p=-\Ls p3,$$ we see that the
result is true for $t\e 1\mod p$. Since $P_{[\f p3]}(-1)=(-1)^{[\f
p3]}P_{[\f p3]}(1)=\sls p3$ and
$$\sum_{x=0}^{p-1}\Ls{x^3-27x+54}p=\sum_{x=0}^{p-1}\Ls{(-3x)^3-27(-3x)+54}p
=\Ls{-3}p\sum_{x=0}^{p-1}\Ls{x^3-3x-2}p=-1,$$ we see that the result
is also true for $t\e -1\mod p$.
\par Now we assume $t\not\e \pm 1\mod p$. Set $W_n(x)=\sum_{k=0}^n\b
nk^2x^k$. From [BM, Theorem 6] we know that
$$W_{[\f p3]}\Big(1-\f x{27}\Big)\e u_p(x)(x-27)^{[\f
p{12}]}J_p\Ls{x(x-24)^3}{x-27}\mod p,$$ where $J_p(x)$ is a certain
Jacobi polynomial given by (3.4) and
$$u_p(x)=\cases 1&\t{if $p\e 1\mod{12}$,}
\\-3(x-24)&\t{if $p\e 5\mod{12}$,}
\\x^2-36x+216&\t{if $p\e 7\mod{12}$,}
\\-3(x-24)(x^2-36x+216)&\t{if $p\e 11\mod{12}$.}
\endcases$$
Set $x=54/(t+1)$. We then have
$$\aligned&W_{[\f p3]}((t-1)/(t+1))
\\&\e \cases \sls {27(1-t)}{1+t}^{[\f
p{12}]}J_p\ls{432(5-4t)^3}{(1-t)(1+t)^3}\mod p&\t{if $p\e
1\mod{12}$,}
\\\f{18(4t-5)}{t+1}\sls {27(1-t)}{1+t}^{[\f
p{12}]}J_p\ls{432(5-4t)^3}{(1-t)(1+t)^3}\mod p&\t{if $p\e
5\mod{12}$,}
\\\f{108(2t^2-14t+11)}{(t+1)^2}\sls {27(1-t)}{1+t}^{[\f
p{12}]}J_p\ls{432(5-4t)^3}{(1-t)(1+t)^3}\mod p&\t{if $p\e
7\mod{12}$,}
\\\f{1944(4t-5)(2t^2-14t+11)}{(t+1)^3}\sls {27(1-t)}{1+t}^{[\f
p{12}]}J_p\ls{432(5-4t)^3}{(1-t)(1+t)^3}\mod p&\t{if $p\e
11\mod{12}$.}\endcases\endaligned\tag 3.5$$ By (3.2) we have
$$W_{[\f p3]}\Ls{t-1}{t+1}
=\Big(1-\f{t-1}{t+1}\Big)^{[\f p3]}P_{[\f
p3]}\Ls{1+(t-1)/(t+1)}{1-(t-1)/(t+1)}=\Ls 2{t+1}^{[\f p3]}P_{[\f
p3]}(t).\tag 3.6$$  If $p\e 2\mod 3$ and $t\e \f 54\mod p$, from the
above we get
$$P_{[\f p3]}\Ls 54=\Ls{\f 54+1}2^{[\f p3]}W_{[\f p3]}\Ls {\f
54-1}{\f 54+1}\e 0\mod p.$$ On the other hand,
$$\sum_{x=0}^{p-1}\Ls{x^3+3(4t-5)x+2(2t^2-14t+11)}p
=\sum_{x=0}^{p-1}\Ls{x^3-27/4}p=\sum_{y=0}^{p-1}\Ls{y-27/4}p=0.$$
Thus the result is true when $p\e 2\mod 3$ and $t\e \f 54\mod p$.
 If $p\e 3\mod 4$
and $2t^2-14t+11\e 0\mod p$, from (3.5) and (3.6) we deduce that
$$P_{[\f p3]}(t)=\Ls{t+1}2^{[\f p3]}W_{[\f p3]}\Ls{t-1}{t+1}
\e 0\mod p.$$ As
$$\sum_{x=0}^{p-1}\Ls{x^3+3(4t-5)x}p=\sum_{x=0}^{p-1}
\Ls{(-x)^3+3(4t-5)(-x)}p=-\sum_{x=0}^{p-1}\Ls{x^3+3(4t-5)x}p,$$ we
see that
$$\sum_{x=0}^{p-1}\Ls{x^3+3(4t-5)x+2(2t^2-14t+11)}p
=\sum_{x=0}^{p-1}\Ls{x^3+3(4t-5)x}p=0.$$ Thus the result is true
when $p\e 3\mod 4$ and $2t^2-14t+11\e 0\mod p$.
 Set $m=3(4t-5)$ and
$n=2(2t^2-14t+11)$. Then
$$4m^3+27n^2=-432(1-t)(1+t)^3\qtq{and so}
\f{2^8\cdot 3^3m^3}{4m^3+27n^2}=\f {432(5-4t)^3}{(1-t)(1+t)^3}.$$ By
the above we need only to assume $m\not\e 0\mod p$ for $p\e 2\mod 3$
and $n\not\e 0\mod p$ for $p\e 3\mod 4$. From (3.3) we see that
$$\align &J_p\Ls{432(5-4t)^3}{(1-t)(1+t)^3}=J_p\Ls{2^8\cdot
3^3m^3}{4m^3+27n^2} \\&\e -(-48m)^{\f{\sls
p3-1}2}(864n)^{\f{\sls{-1}p-1}2}(-16(4m^3+27n^2))^{-[\f p{12}]}
\sum_{x=0}^{p-1}\Ls{x^3+mx+n}p\mod p.\endalign$$
\par If $p\e 1\mod{12}$, from all the above we deduce that
$$\align &P_{[\f p3]}(t)\\&=\Ls{t+1}2^{[\f p3]}W_{[\f p3]}\Ls{t-1}{t+1}
\e \Ls{t+1}2^{\f{p-1}3}\Ls{27(1-t)}{1+t}^{\f
{p-1}{12}}J_p\ls{432(5-4t)^3}{(1-t)(1+t)^3}
\\&\e -2^{-\f{p-1}3}(3(t+1))^{\f{p-1}4}(1-t)^{\f{p-1}{12}}
(16\cdot
432(1-t)(1+t)^3)^{-\f{p-1}{12}}\sum_{x=0}^{p-1}\Ls{x^3+mx+n}p
\\&\e -\sum_{x=0}^{p-1}
\Ls{x^3+3(4t-5)x+2(2t^2-14t+11)}p\mod p.\endalign$$ If $p\e
5\mod{12}$, from all the above we deduce that
$$\align &P_{[\f p3]}(t)\\&=\Ls{t+1}2^{[\f p3]}W_{[\f p3]}\Ls{t-1}{t+1}
\e \Ls{t+1}2^{\f{p-2}3}\f{18(4t-5)}{t+1}\Ls{27(1-t)}{1+t}^{\f
{p-5}{12}}J_p\ls{432(5-4t)^3}{(1-t)(1+t)^3}
\\&\e 2^{-\f{p-5}3}3^{\f{p+3}4}(4t-5)(1+t)^{\f{p-5}4}(1-t)^{\f{p-5}{12}}
 (144(4t-5))^{-1}\\&\qq\times(16\cdot
432(1-t)(1+t)^3)^{-\f{p-5}{12}}\sum_{x=0}^{p-1}\Ls{x^3+mx+n}p
\\&\e \sum_{x=0}^{p-1}
\Ls{x^3+3(4t-5)x+2(2t^2-14t+11)}p\mod p.\endalign$$ If $p\e
7\mod{12}$, from all the above we deduce that
$$\align P_{[\f p3]}(t)&=\Ls{t+1}2^{[\f p3]}W_{[\f p3]}\Ls{t-1}{t+1}
\\&\e\Ls{t+1}2^{\f{p-1}3}\f{108(2t^2-14t+11)}{(t+1)^2}\Ls{27(1-t)}{1+t}^{\f
{p-7}{12}}J_p\ls{432(5-4t)^3}{(1-t)(1+t)^3}
\\&\e -2^{-\f{p-7}3}3^{\f{p+5}4}(2t^2-14t+11)(1+t)^{\f{p-7}4}
(1-t)^{\f{p-7}{12}}(1728(2t^2-14t+11))^{-1}
 \\&\qq\times(16\cdot 432(1-t)(1+t)^3)^{-\f{p-7}{12}}\sum_{x=0}^{p-1}\Ls{x^3+mx+n}p
\\&\e -\sum_{x=0}^{p-1}
\Ls{x^3+3(4t-5)x+2(2t^2-14t+11)}p\mod p.\endalign$$ If $p\e
11\mod{12}$, from all the above we deduce that
$$\align P_{[\f p3]}(t)&=\Ls{t+1}2^{[\f p3]}W_{[\f p3]}\Ls{t-1}{t+1}
\\&\e\Ls{t+1}2^{\f{p-2}3}\f{1944(4t-5)(2t^2-14t+11)}{(t+1)^3}
\Ls{27(1-t)}{1+t}^{\f {p-11}{12}}J_p\ls{432(5-4t)^3}{(1-t)(1+t)^3}
\\&\e 2^{-\f{p-11}3}3^{\f{p-11}4+5}(4t-5)(2t^2-14t+11)(1+t)^{\f{p-11}4}
(1-t)^{\f{p-11}{12}}(48m)^{-1}(864n)^{-1}
 \\&\qq\times(16\cdot 432(1-t)(1+t)^3)^{-\f{p-11}{12}}\sum_{x=0}^{p-1}\Ls{x^3+mx+n}p
\\&\e \sum_{x=0}^{p-1}
\Ls{x^3+3(4t-5)x+2(2t^2-14t+11)}p\mod p.\endalign$$ This proves the
theorem.
\pro{Corollary 3.1} Let $p>3$ be a prime and let $t$ be a
variable. Then
$$\align&\sum_{k=0}^{[p/3]}\b{2k}k\b{3k}k\Ls{1-t}{54}^k
\\&\e P_{[\f p3]}(t)\e -\Ls p3
\sum_{x=0}^{p-1}(x^3+3(4t-5)x+2(2t^2-14t+11))^{\f{p-1}2}\mod
p.\endalign$$
\endpro
Proof. From [S2, Lemma 2.3] we have $P_{[\f p3]}(t)\e
\sum_{k=0}^{[p/3]}\b{2k}k\b{3k}k\sls{1-t}{54}^k\mod p$. By Theorem
3.1 and Euler's criterion, the result is true for
$t=0,1,\ldots,p-1$. Since both sides are polynomials of $t$ with
degree at most $p-1$. Using Lagrange's theorem we obtain the result.
\pro{Corollary 3.2} Let $p\ge 17$ be a prime and $t\in\Bbb Z_p$.
Then
$$\align&\sum_{x=0}^{p-1}\Ls{x^3+3(4t-5)x+2(2t^2-14t+11)}p
\\&=\Ls p3 \sum_{x=0}^{p-1}\Ls{x^3-3(4t+5)x+2(2t^2+14t+11)}p.\endalign$$
\endpro
Proof. Since $P_{[\f p3]}(-t)=(-1)^{[\f p3]}P_{[\f p3]}(t)=\ls p3
P_{[\f p3]}(t)$, by Theorem 3.1 we have
$$\align&\sum_{x=0}^{p-1}\Ls{x^3+3(4t-5)x+2(2t^2-14t+11)}p
\\&\e \Ls p3
\sum_{x=0}^{p-1}\Ls{x^3-3(4t+5)x+2(2t^2+14t+11)}p\mod p.\endalign$$
By Weil's estimate ([BEW, p.183]) we have
$$\align&\Big|\sum_{x=0}^{p-1}\Ls{x^3+3(4t-5)x+2(2t^2-14t+11)}p\Big|\le
2\sqrt
p,\\&\Big|\sum_{x=0}^{p-1}\Ls{x^3-3(4t+5)x+2(2t^2+14t+11)}p\Big|\le
2\sqrt p.\endalign$$ Since $4\sqrt p<p$ for $p\ge 17$, from the
above we deduce the result.

\pro{Corollary 3.3} Let $p>3$ be a prime. Then
$$\sum_{x=0}^{p-1}\Ls{x^3-120x+506}p=\cases \sls 2pL&\t{if $3\mid
p-1$, $4p=L^2+27M^2$ and $3\mid L-1$,}\\0&\t{if $p\e 2\mod 3$.}
\endcases$$\endpro
Proof. It is easy to check the result for $p=5,7,11,13$. Now we
assume $p\ge 17$. Taking $t=\f 54$ in Corollary 3.2 we find that
$$\align\sum_{x=0}^{p-1}\Ls{x^3-\f{27}4}p&=\Ls p3
\sum_{x=0}^{p-1}\Ls{x^3-30x+\f{253}4}p =\Ls
p3\sum_{x=0}^{p-1}\Ls{(\f x2)^3-30\cdot \f x2+\f{253}4}p
\\&=\Ls p3\Ls 2p\sum_{x=0}^{p-1}\Ls{x^3-120x+506}p.\endalign$$
For $p\e 2\mod 3$ it is clear that
$$\sum_{x=0}^{p-1}\Ls{x^3-\f{27}4}p=\sum_{x=0}^{p-1}\Ls{x-\f{27}4}p
=\sum_{x=0}^{p-1}\Ls xp=0.$$ Thus the result is true when $p\e 2\mod
3$.
\par Now assume $p\e 1\mod 3$, $p=A^2+3B^2$, $4p=L^2+27M^2$ and $A\e
L\e 1\mod 3$. It is known that $2^{\f{p-1}3} \e 1\mod p$ if and only
if $3\mid B$. When $3\nmid B$ we choose the sign of $B$ so that $B\e
1\mod 3$. By [S1, (2.12)] we have $2^{(p-1)/3}\e \f 12(-1-\f AB)\mod
p$. From [S1, (2.9)-(2.11)] we deduce that
$$\aligned  \sum_{x=0}^{p-1}\Ls{x^3-27/4}p
&=1+\sum_{x=1}^{p-1}\Ls{x^3-27/4}p
\\&=\cases -2A=L\mod p&\t{if $2^{\f{p-1}3}\e 1\mod p$,}
\\A+3B=L\mod p&\t{if $2^{\f{p-1}3}\not\e 1\mod p$ and $B\e 1\mod 3$.}
\endcases\endaligned$$
Thus
$$\sum_{x=0}^{p-1}\Ls{x^3-120x+506}p=\Ls 2p\sum_{x=0}^{p-1}\Ls{x^3-27/4}p
=\Ls 2pL.$$ This completes the proof.

\pro{Theorem 3.2} Let $p>3$ be a prime. Then
\par $(\t{\rm i})$ If $p\e 2\mod 3$, then
$$\sum_{k=0}^{[p/3]}\f{\b{2k}k\b{3k}k}{(-216)^k}
\e \sum_{k=0}^{[p/3]}\f{\b{2k}k\b{3k}k}{24^k}\e P_{[\f p3]}\Ls 54 \e
0\mod p.$$
\par $(\t{\rm ii})$ If $p\e 1\mod 3$ and so $4p=L^2+27M^2$
with $L,M\in\Bbb Z$ and $L\e 1\mod 3$, then
$$ \sum_{k=0}^{[p/3]}\f{\b{2k}k\b{3k}k}{(-216)^k}
\e\sum_{k=0}^{[p/3]}\f{\b{2k}k\b{3k}k}{24^k} \e P_{[\f p3]}\Ls 54 \e
-L \e \Ls{-2}p\b{\f{2(p-1)}3}{[\f p{12}]}\mod p.$$
\endpro
Proof. Putting $t=\pm\f 54$ in Corollary 3.1 we get
$$P_{[\f p3]}\Ls 54\e \sum_{k=0}^{[p/3]}\f{\b{2k}k\b{3k}k}{(-216)^k}
\mod p\qtq{and}P_{[\f p3]}\Big(-\f 54\Big)\e
\sum_{k=0}^{[p/3]}\f{\b{2k}k\b{3k}k}{24^k}\mod p.$$ This together
with (1.3) yields
$$\sum_{k=0}^{[p/3]}\f{\b{2k}k\b{3k}k}{(-216)^k} \e\Ls
p3\sum_{k=0}^{[p/3]}\f{\b{2k}k\b{3k}k}{24^k}\mod p.$$ If $p\e 2\mod
3$, by Theorem 3.1 we have
$$P_{[\f p3]}\Ls 54\e \sum_{x=0}^{p-1}\Ls{x^3-\f{27}4}p
=\sum_{x=0}^{p-1}\Ls{x-\f{27}4}p =\sum_{x=0}^{p-1}\Ls xp=0 \mod p.$$
Thus (i) is true.
\par Now assume $p\e 1\mod 3$, $4p=L^2+27M^2$ and $L\e 1\mod 3$.
By Theorem 3.1 and the proof of Corollary 3.3 we have
$$P_{[\f p3]}\Ls 54\e -\sum_{x=0}^{p-1}\Ls{x^3-\f{27}4}p=-L\mod p.$$
On the other hand, by the proof of Theorem 3.1,
$$\aligned &P_{[\f p3]}\Ls 54=\Ls{\f 54+1}2^{[\f p3]}W_{[\f p3]}\Ls {\f
54-1}{\f 54+1}\\&\e\cases \ls 98^{\f{p-1}3}\ls{27(1-\f 54)}{1+\f
54}^{\f{p-1}{12}}J_p(0)\e (-1)^{\f{p-1}{12}}3^{-\f{p-1}4}J_p(0)\mod
p\q\t{if $p\e 1\mod{12}$,}
\\\ls 98^{\f{p-1}3}\f{108(2\sls 54^2-14\cdot \f 54+11)}{(\f 54+1)^2}
\ls{27(1-\f 54)}{1+\f 54}^{\f{p-7}{12}}J_p(0)
\\\q\e
-8(-1)^{\f{p-7}{12}}3^{-\f{p-7}4}J_p(0)\mod p\qq\qq\qq\qq\qq\t{if
$p\e 7\mod{12}$.}
\endcases\endaligned$$
By the definition of $J_p(x)$, we have
$$\align J_p(0)&=1728^{[\f p{12}]}\cdot 2^{-[\f p{12}]}\sum_{r=0}^{[\f
p{12}]}\b{[\f p{12}]-\f 13\sls p3}r\b{[\f p{12}]-\f 12\sls{-1}p}{[\f
p{12}]-r}0^{[\f p{12}]-r}2^r \\&=1728^{[\f p{12}]}\b{[\f p{12}]-\f
13\ls p3}{[\f p{12}]}=(-1728)^{[\f p{12}]}\b{\f 13\sls p3-1}{[\f
p{12}]}.\endalign$$ Hence
$$J_p(0)\e (-1728)^{[\f p{12}]}\b{\f{2(p-1)}3}{[\f
p{12}]}\mod p$$ and therefore
$$\aligned P_{[\f p3]}\Ls 54\e\cases
(-1)^{\f{p-1}{12}}3^{-\f{p-1}4}(-1728)^{\f{p-1}{12}}
\b{\f{2(p-1)}3}{\f{p-1}{12}}\e \ls
2p\b{\f{2(p-1)}3}{\f{p-1}{12}}\mod p &\t{if $12\mid p-1$,}
\\-8(-1)^{\f{p-7}{12}}3^{-\f{p-7}4}(-1728)^{\f{p-7}{12}}
\b{\f{2(p-1)}3}{\f{p-7}{12}} \e -\ls 2p\b{\f{2(p-1)}3}{\f{p-7}{12}}
\mod p &\t{if $12\mid p-7$.}\endcases\endaligned$$ \par Now putting
all the above together we obtain the result.
\par\q
\newline{\bf Remark 3.1} For any prime $p>3$, Zhi-Wei Sun
conjectured that ([Su1, Conjecture A46])
$$\sum_{k=0}^{p-1}\f{(3k)!}{24^k\cdot k!^3} \e\Ls p3
\sum_{k=0}^{p-1}\f{(3k)!}{(-216)^k\cdot k!^3} \e \cases
\b{2(p-1)/3}{(p-1)/3}\mod{p^2}&\t{if $p\e 1\mod 3$,}
\\0\mod p&\t{if $p\e 2\mod 3$.}
\endcases$$

 \subheading{4. Congruences for $\sum_{k=0}^{p-1}
 \f{\b{2k}k^2\b{3k}k}{m^k}$}
 \par Let $p>3$ be a prime and $m\in\Bbb Z$ with $p\nmid m$.
  In the section we partially solve Z.W.
Sun's conjectures on $\sum_{k=0}^{p-1}\b{2k}k^2 \b{3k}km^{-k}\mod
{p^2}$.
 \pro{Theorem 4.1} Let $p>3$ be
a prime, $m\in\Bbb Z_p$, $m\not\e 0\mod p$ and $t=\sqrt{1-108/m}$.
Then
$$\sum_{k=0}^{p-1}\f{\b{2k}k^2\b{3k}k}{m^k}\e
P_{[\f p3]}(t)^2\e
\Big(\sum_{x=0}^{p-1}(x^3+3(4t-5)x+2(2t^2-14t+11))^{\f{p-1}2}\Big)^2\mod
p.$$ Moreover, if $P_{[\f p3]}(t)\e 0\mod p$ or
$\sum_{x=0}^{p-1}(x^3+3(4t-5)x+2(2t^2-14t+11))^{\f{p-1}2}\e 0\mod
p$, then
$$\sum_{k=0}^{p-1}\f{\b{2k}k^2\b{3k}k}{m^k}\e 0\mod{p^2}.$$
\endpro
Proof. Since $\f{1-t}{54}(1-27\cdot \f{1-t}{54})=\f 1m$, by Theorem
2.1 we have
$$\sum_{k=0}^{p-1}\f{\b{2k}k^2\b{3k}k}{m^k}\e
\Big(\sum_{k=0}^{p-1}\b{2k}k\b{3k}k\Ls{1-t}{54}^k\Big)^2\mod{p^2}.
\tag 4.1$$ Observe that $p\mid \b{2k}k\b{3k}k$ for $[\f p3]<k<p$.
From the above and Corollary 3.1 we see that
$$\align&\sum_{k=0}^{p-1}\b{2k}k\b{3k}k\Ls{1-t}{54}^k
\\&\e P_{[\f p3]}(t)\e -\Ls
p3\sum_{x=0}^{p-1}(x^3+3(4t-5)x+2(2t^2-14t+11))^{\f{p-1}2}\mod
p.\endalign$$ This together with (4.1) yields the result.

 \pro{Theorem 4.2 ([Su1, Conjecture A8])} Let $p>3$ be a
prime. Then
$$\sum_{k=0}^{p-1}\f{\b{2k}k^2\b{3k}k}{(-192)^k}
\e\cases L^2\mod p&\t{if $p\e 1\mod 3$ and so $4p=L^2+27M^2$,}
\\0\mod {p^2}&\t{if $p\e 2\mod 3$}.
\endcases$$
\endpro
Proof. Putting $m=-192$ and $t=\f 54$ in Theorem 4.1 and then
applying Theorem 3.2 we obtain the result.

\pro{Lemma 4.1} Let $p$ be an odd prime and let $a,m,n$ be p-adic
integers. Then
$$\sum_{x=0}^{p-1}(x^3+a^2mx+a^3n)^{\f{p-1}2}\e a^{\f{p-1}2}
\sum_{x=0}^{p-1}(x^3+mx+n)^{\f{p-1}2}\mod p.$$ Moreover, if $a,m,n$
are congruent to some integers, then
$$\sum_{x=0}^{p-1}\Ls{x^3+a^2mx+a^3n}p=
\Ls ap\sum_{x=0}^{p-1}\Ls{x^3+mx+n}p.$$
\endpro
Proof. For any positive integer $k$ it is well known that (see [IR])
$$\sum_{x=0}^{p-1}x^k\e \cases p-1\mod p&\t{if $p-1\mid k$,}
\\0\mod p&\t{if $p-1\nmid k$.}
\endcases$$ Since
$$\align&\sum_{x=0}^{p-1}(x^3+a^2mx+a^3n)^{\f{p-1}2}
\\&=\sum_{x=0}^{p-1}\sum_{k=0}^{(p-1)/2}\b{(p-1)/2}k
(x^3+a^2mx)^k(a^3n)^{\f{p-1}2-k}
\\&=\sum_{x=0}^{p-1}\sum_{k=0}^{(p-1)/2}\b{(p-1)/2}k
\sum_{r=0}^k\b krx^{3r}(a^2mx)^{k-r}(a^3n)^{\f{p-1}2-k}
\\&=\sum_{r=0}^{(p-1)/2}\sum_{k=r}^{(p-1)/2}
\b{(p-1)/2}k\b
kr(a^2m)^{k-r}(a^3n)^{\f{p-1}2-k}\sum_{x=0}^{p-1}x^{k+2r}
\\&\e (p-1)\sum_{r=0}^{(p-1)/2}\b{(p-1)/2}{p-1-2r}\b{p-1-2r}
r(a^2m)^{p-1-3r}(a^3n)^{2r-\f{p-1}2}
\\&=a^{\f{p-1}2}(p-1)\sum_{\f{p-1}4\le r\le \f{p-1}3}\b{(p-1)/2}{p-1-2r}\b{p-1-2r}
rm^{p-1-3r}n^{2r-\f{p-1}2}\mod p,
\endalign$$
we see that the congruence in Lemma 4.1 is true.
\par Now suppose that $a,m,n$
are congruent to some integers. If $a\e 0\mod p$, clearly
$$\sum_{x=0}^{p-1}\Ls{x^3+a^2mx+a^3n}p=\sum_{x=0}^{p-1}\Ls{x^3}p
=\sum_{x=0}^{p-1}\Ls xp=0=\Ls ap\sum_{x=0}^{p-1}\Ls{x^3+mx+n}p.$$ If
$a\not\e 0\mod p$, then clearly
$$\sum_{x=0}^{p-1}\Ls{x^3+a^2mx+a^3n}p=\sum_{x=0}^{p-1}\Ls{(ax)^3+a^2m(ax)+a^3n}p
=\Ls ap\sum_{x=0}^{p-1}\Ls{x^3+mx+n}p.$$ Thus the lemma is proved.

\pro{Lemma 4.2} Let $p$ be an odd prime. Then
$$\aligned &\sum_{x=0}^{p-1}\Ls{x^3-30x-56}p
\\&= \cases (-1)^{[\f p8]+1}\sls 3p2c&\t{if $p\e 1,3\mod 8$,
$p=c^2+2d^2$ and $4\mid c-1$,}
\\0&\t{if $p\e 5,7\mod 8$.}
\endcases\endaligned$$
\endpro
Proof. From [BE, Theorems 5.12 and 5.17] we know that
$$\sum_{k=0}^{p-1}\Ls{x^3-4x^2+2x}p=\cases (-1)^{[\f p8]+1}2c&\t{if $p=c^2+2d^2\e
1,3\mod 8$ with $4\mid c-1$,}
\\0&\t{otherwise.}\endcases$$
As $27(x^3-4x^2+2x)=(3x-4)^3-30(3x-4)-56$, we see that
$$\align&\sum_{x=0}^{p-1}\Ls{x^3-4x^2+2x}p
\\&=\Ls 3p\sum_{x=0}^{p-1}\Ls{(3x-4)^3-30(3x-4)-56}p =\Ls
3p\sum_{x=0}^{p-1}\Ls{x^3-30x-56}p .
\endalign$$ Thus the result follows.

\pro{Lemma 4.3} Let $p$ be an odd prime. Then
$$\aligned&\sum_{n=0}^{p-1}(n^3-(15+30\sqrt{-2})n-28+70\sqrt{-2}))^{\f{p-1}2}
\\&\e \cases \sls{2+\sqrt{-2}}p(-1)^{[\f p8]+1}\sls 3p2c
\mod p&\t{if $p=c^2+2d^2\e 1,3\mod 8$ and $4\mid c-1$,}
\\0\mod p&\t{if $p\e 5,7\mod 8$.}\endcases\endaligned$$
\endpro
Proof. It is easily seen that
$$-15(1+2\sqrt {-2})=-30\Ls{1-\sqrt{-2}}{\sqrt{-2}}^2\qtq{and}
-28+70\sqrt{-2}=-56 \Ls{1-\sqrt{-2}}{\sqrt{-2}}^3.$$ Thus, by Lemmas
4.1 and 4.2 we have
$$\aligned&\sum_{n=0}^{p-1}\big(n^3-(15+30\sqrt{-2})n-28+70\sqrt{-2})
\big)^{\f{p-1}2}
\\&\e\Ls{1-\sqrt{-2}}{\sqrt{-2}}^{\f{p-1}2}
\sum_{n=0}^{p-1}(n^3-30n-56)^{\f{p-1}2}
\\&\e \Big(-\f{2+\sqrt{-2}}2\Big)^{\f{p-1}2}
\sum_{n=0}^{p-1}\Ls{n^3-30n-56}p \\&
\e\cases\sls{2+\sqrt{-2}}p(-1)^{[\f p8]+1}\sls 3p2c\mod p&\t{if
$p=c^2+2d^2\e 1,3\mod 8$ and $4\mid c-1$,}
\\0\mod p&\t{if $p\e 5,7\mod 8$.}\endcases\endaligned$$
This proves the lemma.
 \pro{Theorem
4.3} Let $p$ be an odd prime. Then
$$\aligned&P_{[\f p3]}(5/\sqrt {-2})\\&\e\cases
(-1)^{[\f{p}8]}\sls{-2-\sqrt{-2}}p2c \mod p&\t{if $p=c^2+2d^2\e
1,3\mod 8$ and $4\mid c-1$,}\\0\mod p&\t{if $p\e 5,7\mod
8$}\endcases\endaligned$$ and
$$\sum_{k=0}^{p-1}\f{\b{2k}k^2\b{3k}k}{8^k}\e
 \cases 4c^2\mod p&\t{if $p=c^2+2d^2\e 1,3\mod 8$,}
 \\0\mod{p^2}&\t{if $p\e 5,7\mod 8$.}\endcases$$
\endpro
Proof. From Corollary 3.1 and Lemma 4.3 we deduce that
$$\aligned P_{[\f p3]}\Ls{5}{\sqrt {-2}}&\e
-\Ls
p3\sum_{n=0}^{p-1}(n^3+3(-10\sqrt{-2}-5)n-28+70\sqrt{-2})^{\f{p-1}2}
\\&
\e \cases \sls p3\sls{2+\sqrt{-2}}p(-1)^{[\f p8]}\sls 3p2c
=\sls{-2-\sqrt{-2}}p(-1)^{[\f p8]}2c\mod p\\\qq\qq\qq\; \t{if
$p=c^2+2d^2\e 1,3\mod 8$ and $4\mid c-1$,}
\\0\mod p\qq\t{if $p\e 5,7\mod 8$.}\endcases\endaligned$$
 Now taking $m=8$ and $t=5/\sqrt{-2}$ in
Theorem 4.1 and then applying the above we deduce the remaining
result.
\newline{\bf Remark 4.1} Let $p$ be an odd prime. Zhi-Wei Sun
 conjectured that ([Su1,
 Conjecture A5])
$$\sum_{k=0}^{p-1}\f{\b{2k}k^2\b{3k}k}{8^k}\e \cases 4c^2-2p\mod{p^2}
&\t{if $p=c^2+2d^2\e 1,3\mod 8$,}\\0\mod{p^2}&\t{if $p\e 5,7\mod
8$.}\endcases$$

\pro{Corollary 4.1} Let $p$ be a prime such that $p\e 1,3\mod 8$ and
$p=c^2+2d^2$ with $c\e 1\mod 4$. Then
$$P_{[\f p3]}\Ls {5d}c
\e\cases \sls{c-d}3 2c\mod p&\t{if $p\e 1\mod 8$,}
\\(-1)^{\f{d-1}2}\sls{c-d}32c\mod p&\t{if $p\e 3\mod 8$.}
\endcases$$\endpro
Proof. It is clear that $\sls dp=\sls pd=\sls{c^2}d=1$. Thus
$$\aligned
\Ls{2+c/d}p&=\Ls dp\sls{c+2d}p=\Ls{c+2d}p=(-1)^{\f{p-1}2\cdot
\f{2d+c-1}2}\Ls p{c+2d}\\&=(-1)^{\f{p-1}2d}\Ls
{(c+2d)(c-2d)+6d^2}{c+2d} =(-1)^{\f{p-1}2d}\Ls
{6}{c+2d}\\&=(-1)^{\f{p-1}2d+\f{(c+2d)^2-1}8}\Ls{3}{c+2d}
=(-1)^{\f{p-1+4d+2d^2}8}\cdot
(-1)^{\f{p-1}2d+\f{c+2d-1}2}\Ls{c+2d}3\\&=
(-1)^{\f{p-1+4d+2d^2}8}\Ls{c-d}3 =\cases
(-1)^{\f{p-1}8}\sls{c-d}3&\t{if $p\e 1\mod 8$,}
\\(-1)^{\f{p-3}8+\f{d+1}2}\sls{c-d}3&\t{if $p\e 3\mod 8$.}
\endcases
\endaligned$$
By Theorem 4.3 we have
$$P_{[\f p3]}\Ls 5{c/d}\e (-1)^{[\f p8]}\Ls{-1}p\Ls{2+c/d}p2c\mod
p.$$ Now combining all the above we obtain the result.

\pro{Lemma 4.4} Let $p$ be an odd prime and $p\not=11$. Then
$$\sum_{x=0}^{p-1}\Ls{x^3-24\cdot 11x+14\cdot 11^2}p
=\cases \sls u{11}u&\t{if $\sls p{11}=1$ and so
$4p=u^2+11v^2$,}\\0&\t{if $\sls p{11}=-1$.}\endcases$$
\endpro
Proof.  It is known that (see [PR] and [JM])
$$\sum_{x=0}^{p-1}\Ls{x^3-96\cdot 11x+112\cdot 11^2}p=\cases \sls 2p\sls
u{11}u&\t{if $\sls p{11}=1$ and $4p=u^2+11v^2$,}\\0&\t{if $\sls
p{11}=-1$.}\endcases$$ Since
$$\align&\sum_{x=0}^{p-1}\Ls{x^3-96\cdot 11x+112\cdot 11^2}p
\\&=\sum_{x=0}^{p-1}\Ls{(2x)^3-96\cdot 11\cdot 2x+112\cdot 11^2}p
 =\Ls
2p\sum_{x=0}^{p-1}\Ls{x^3-24\cdot 11x+14\cdot 11^2}p,\endalign$$ we
deduce the result.

\pro{Lemma 4.5} Let $p\not=11$ be an odd prime. Then
$$\aligned&\sum_{n=0}^{p-1}(n^3+12(-5+\sqrt{-11})n
+14(11-4\sqrt{-11}))^{\f{p-1}2} \\&\e\cases
\sls{-22+2\sqrt{-11}}p\sls u{11}u\mod p&\t{if $\sls p{11}=1$ and so
$4p=u^2+11v^2$,}\\0\mod p&\t{if $\sls
p{11}=-1$.}\endcases\endaligned$$
\endpro
Proof. It is easily seen that
$$12(-5+\sqrt{-11})=-24\cdot
11\Ls{\sqrt{-11}+1}{2\sqrt{-11}}^2\ \t{and}\
14(11-4\sqrt{-11})=14\cdot 11^2 \Ls{\sqrt{-11}+1}{2\sqrt{-11}}^3.$$
Thus, by Lemma 4.1 we have
$$\align&\sum_{n=0}^{p-1}\big(n^3+12(-5+\sqrt{-11})n+14(11-4\sqrt{-11})
\big)^{\f{p-1}2}
\\&\e \Ls{\sqrt{-11}+1}{2\sqrt{-11}}^{\f{p-1}2}
\sum_{x=0}^{p-1}\big(x^3-24\cdot 11x+14\cdot 11^2\big)^{\f{p-1}2}
\\&\e\Ls{-22+2\sqrt{-11}}{-11\cdot 4}^{\f{p-1}2}\sum_{x=0}^{p-1}\Ls{x^3-24\cdot 11x+14\cdot
11^2}p\mod p.\endalign$$ Now applying Lemma 4.4 we deduce the
result.

  \pro{Theorem 4.4} Let
$p\not=11$ be an odd prime. Then
$$P_{[\f p3]}\Ls {\sqrt {-11}}4
\e \cases -\sls p3 \sls{-11+\sqrt{-11}}p\sls u{11}u\mod p&\t{if
$\sls p{11}=1$ and so $4p=u^2+11v^2$,}\\0\mod p&\t{if $\sls
p{11}=-1$}
\endcases$$ and
$$\sum_{k=0}^{p-1}\f{\b{2k}k^2\b{3k}k}{64^k}\e\cases
 u^2\mod p&\t{if $\sls
p{11}=1$ and so $4p=u^2+11v^2$,}\\0\mod {p^2}&\t{if $\sls
p{11}=-1$.}\endcases$$
\endpro
Proof. From Corollary 3.1 and Lemma 4.5 we deduce that
$$\aligned P_{[\f p3]}\Ls{\sqrt {-11}}4&\e
-\Ls p3\sum_{n=0}^{p-1}\Big(n^3+3(\sqrt{-11}-5)n+\f{-11}4+22-7
\sqrt{-11}\Big)^{\f{p-1}2}
\\&\e -\Ls
p3\sum_{n=0}^{p-1}\Big(\ls n2^3+3(\sqrt{-11}-5)\f
n2+\f{77-28\sqrt{-11}}4\Big)^{\f{p-1}2}
\\&\e -\Ls p3\Ls 2p\sum_{n=0}^{p-1}\big(n^3+12(-5+\sqrt{-11})n+14(11-4
\sqrt{-11})\big)^{\f{p-1}2} \\&\e\cases -\sls p3
\sls{-11+\sqrt{-11}}p\sls u{11}u\mod p&\t{if $\sls p{11}=1$ and so
$4p=u^2+11v^2$,}\\0\mod p&\t{if $\sls
p{11}=-1$.}\endcases\endaligned$$ Now taking $m=64$ and $t=\f{\sqrt
{-11}}4$ in Theorem 4.1 and then applying the above we deduce the
remaining result.
\newline{\bf Remark 4.2} Let $p$ be an odd prime such that $p\not =11$.
 Zhi-Wei Sun conjectured
that ([Su1, Conjecture A4])
$$\sum_{k=0}^{p-1}\f{\b{2k}k^2\b{3k}k}{64^k}\e \cases u^2-2p\mod{p^2}
&\t{if $\sls p{11}=1$ and so $4p=u^2+11v^2$,} \\0\mod{p^2}&\t{if
$\sls p{11}=-1$.}\endcases$$

\pro{Corollary 4.2} Let $p$ be an odd prime such that $\sls p{11}=1$
and so $4p=u^2+11v^2$ with $u,v\in\Bbb Z$. Let $v=2^{\alpha}v_0$ and
$u-11v=2^{\beta}w$ with $2\nmid v_0w$. Then
$$P_{[\f p3]}\Ls u{4v}\e (-1)^{\f{p-1}2\cdot \f{v_0-w}2}
\Ls 2p^{\alpha+\beta}\Ls{v-u}3u\mod p.$$
\endpro
Proof. It is clear that
$$\align \Ls{-11+u/v}p&=\Ls{v(u-11v)}p=\Ls 2p^{\alpha+\beta}\Ls{v_0}p\Ls wp
\\&=\Ls 2p^{\alpha+\beta}(-1)^{\f{v_0-1}2\cdot \f{p-1}2}\Ls p{v_0}
\cdot (-1)^{\f{w-1}2\cdot\f{p-1}2}\Ls pw
\\&=\Ls 2p^{\alpha+\beta}(-1)^{\f{v_0-w}2\cdot \f{p-1}2}\Ls {u^2+11v^2}
{v_0}\Ls{u^2+11v^2}w \\&=\Ls 2p^{\alpha+\beta}(-1)^{\f{v_0-w}2\cdot
\f{p-1}2}\Ls {u^2} {v_0}\Ls{u^2-121v^2+132v^2}w
\\&=\Ls 2p^{\alpha+\beta}(-1)^{\f{v_0-w}2\cdot
\f{p-1}2}\Ls{33}w =\Ls 2p^{\alpha+\beta}(-1)^{\f{v_0-w}2\cdot
\f{p-1}2}\Ls {u-11v}{33}
\\&=\Ls 2p^{\alpha+\beta}(-1)^{\f{v_0-w}2\cdot
\f{p-1}2}\Ls{u+v}3\Ls u{11}.\endalign$$ Thus, by Theorem 4.4 and the
above we get
$$\align P_{[\f p3]}\Ls u{4v}&\e -\Ls p3\Ls{-11+u/v}p\Ls u{11}u
\e -\Ls{u^2+11v^2}3\Ls 2p^{\alpha+\beta}(-1)^{\f{v_0-w}2\cdot
\f{p-1}2}\Ls{u+v}3u
\\&=\Ls 2p^{\alpha+\beta}(-1)^{\f{v_0-w}2\cdot
\f{p-1}2}\Ls{v-u}3u\mod p.\endalign$$ This proves the corollary.

\par\q
\par Let $p>3$ be a prime and let $\Bbb F_p$ be the field of $p$ elements.
For $m,n\in \Bbb F_p$ let $\#E_p(x^3+mx+n)$ be the number of points
on the curve $E$: $y^2=x^3+mx+n$ over the field $\Bbb F_p$. It is
well known that (see for example [S1, pp.221-222])
$$\#E_p(x^3+mx+n)=p+1+\sum_{x=0}^{p-1}\Ls{x^3+mx+n}p.\tag 4.2$$
Let $K=\Bbb Q(\sqrt{-d})$ be an imaginary quadratic field and the
curve $y^2=x^3+mx+n$ has complex multiplication by $K$. By Deuring's
theorem ([C, Theorem 14.16],[PV],[I]), we have
$$\#E_p(x^3+mx+n)=\cases p+1&\t{if $p$ is inert in $K$,}
\\p+1-\pi-\bar{\pi}&\t{if $p=\pi\bar{\pi}$ in $K$,}
\endcases\tag 4.3$$ where $\pi$ is in an order in $K$ and
$\bar{\pi}$ is the conjugate number of $\pi$. If $4p=u^2+dv^2$ with
$u,v\in\Bbb Z$, we may take $\pi=\f 12(u+v\sqrt{-d})$. Thus,
$$\sum_{x=0}^{p-1}\Ls{x^3+mx+n}p=\cases \pm u&\t{if $4p=u^2+dv^2$
with $u,v\in\Bbb Z$,}\\0&\t{otherwise.}
\endcases\tag 4.4$$
In [Gr], [JM] and [PV] the sign of $u$ in (4.4) was determined for
those imaginary quadratic fields $K$ with class number $1$. In [LM]
and [I] the sign of $u$ in (4.4) was determined for imaginary
quadratic fields $K$ with class number $2$. For general results on
the sign of $u$ in (4.4), see [M], [St], [RS] and the survey [Ri].

\pro{Lemma 4.6} Let $p$ be a prime with $p\e \pm 1\mod 8$. Then
$$\aligned&\sum_{n=0}^{p-1}\Ls{n^3+(-15+6\sqrt 2)n+24-14\sqrt 2}p
\\&=\cases 2x\sls {2x}3&\t{if $p\e 1,7\mod {24}$ and so $p=x^2+6y^2$,}
\\0&\t{if $p\e 17,23\mod{24}.$}\endcases\endaligned$$
\endpro
Proof. It is easy to check the result for $p=7$. Now assume $p\ge
17$. From [I, p.133] we know that the elliptic curve defined by the
equation $y^2=x^3+(-21+12\sqrt 2)x-28+22\sqrt 2$ has complex
multiplication by the order of discriminant $-24$. Since
$4p=u^2+24v^2$ implies $2\mid u$ and $p=(\f u2)^2+6v^2$, by (4.4)
and [I, Theorem 3.1] we have
$$\aligned&\sum_{n=0}^{p-1}\Ls{n^3+(-21+12\sqrt 2)n-28+22\sqrt 2}p
\\&=\cases 2x\sls {2x}3\sls{1+\sqrt 2}p&\t{if $p\e 1,7\mod {24}$
 and so $p=x^2+6y^2$,}
\\0&\t{if $p\e 17,23\mod{24}.$}\endcases\endaligned$$
Observe that
$$\f{-15-6\sqrt 2}{-21+12\sqrt 2}=(1+\sqrt 2)^2\qtq{and}
\f{24+14\sqrt 2}{-28+22\sqrt 2}=(1+\sqrt 2)^3.$$ Using Corollary 3.2
and Lemma 4.1 we see that
$$\align&\Ls p3\sum_{n=0}^{p-1}\Ls{n^3+(-15+6\sqrt 2)n+24-14\sqrt 2}p
\\&=\sum_{n=0}^{p-1}\Ls{n^3-(15+6\sqrt 2)n+24+14\sqrt 2}p
\\&=\Ls{1+\sqrt 2}p\sum_{n=0}^{p-1}\Ls{n^3+
(-21+12\sqrt 2)n-28+22\sqrt 2}p.\endalign$$
  Now putting all the above together we
obtain the result.

\pro{Theorem 4.5} Let $p$ be a prime such that $p\e 1,7\mod 8$. Then
$$P_{[\f p3]}\Ls{\sqrt 2}2\e \cases
2x\sls x3\mod p&\t{if $p=x^2+6y^2\e 1,7\mod {24}$,}
 \\0\mod p&\t{if
$p\e 17,23\mod{24}$}\endcases$$ and
$$\sum_{k=0}^{p-1}\f{\b{2k}k^2\b{3k}k}{216^k}\e
\cases 4x^2\mod p&\t{if $p=x^2+6y^2\e 1,7\mod {24}$,}
\\0\mod {p^2}&\t{if $p\e 17,23\mod {24}$.}
\endcases$$
\endpro
Proof. From Theorem 3.1 and Lemma 4.6 we deduce that
$$\aligned P_{[\f p3]}\Ls{\sqrt 2}2&\e
-\Ls p3\sum_{n=0}^{p-1}\Ls{n^3+(-15+6\sqrt 2)n+24-14\sqrt 2}p
\\&\e \cases -2x\sls{2x}3=2x\sls x3\mod p&\t{if $p=x^2+6y^2\e
1,7\mod{24}$,}
\\0\mod p&\t{if $p\e 17,23\mod{24}$.}\endcases
\endaligned$$
Taking $m=216$ and $t=\f{\sqrt 2}2$ in Theorem 4.1 and then applying
the above we deduce the remaining result.
\newline{\bf Remark 4.3} For any prime $p>3$, Z.W. Sun conjectured
that ([Su1, A14])
$$\sum_{k=0}^{p-1}\f{\b{2k}k^2\b{3k}k}{216^k}\e
\cases 4x^2-2p\mod {p^2}&\t{if $p=x^2+6y^2\e 1,7\mod {24}$,}
\\8x^2-2p\mod {p^2}&\t{if $p=2x^2+3y^2\e 5,11\mod {24}$,}
\\0\mod {p^2}&\t{if $p\e 13,17,19,23\mod {24}$.}
\endcases$$
 \pro{Conjecture 4.1} Let
$p$ be a prime such that $p\e 5,11,13,19\mod{24}$. Then
$$P_{[\f p3]}\Ls {\sqrt 2}2\e \cases 0\mod p&\t{if $p\e
13,19\mod{24}$,}
\\-2x\sls x3\sqrt 2\mod p&\t{if $p=2x^2+3y^2\e 5,11\mod {24}$.}
\endcases$$
\endpro
\pro{Lemma 4.7} Let $p$ be a prime with $p\e \pm 1\mod 5$. Then
$$\aligned&\sum_{n=0}^{p-1}\Ls{n^3+(-15+12\sqrt 5)n+42-28\sqrt 5}p
\\&=\cases 2x\sls {2x}3&\t{if $p\e 1,4\mod {15}$ and so $p=x^2+15y^2$,}
\\0&\t{if $p\e 11,14\mod{15}.$}\endcases\endaligned$$
\endpro
Proof. From [I, Proposition 3.3] we know that the elliptic curve
defined by the equation $y^2=x^3+(105+48\sqrt 5)x-784-350\sqrt 5$
has complex multiplication by the order of discriminant $-15$. Since
$4p=u^2+60v^2$ implies $2\mid u$ and $p=(\f u2)^2+15v^2$, by (4.4)
and [I, Theorem 3.1] we have
$$\aligned&\sum_{n=0}^{p-1}\Ls{n^3+(105+48\sqrt 5)n-784-350\sqrt 5}p
\\&=\cases 2x\sls {2x}3\sls{(1+\sqrt 5)/2}p&\t{if $p\e 1,4\mod {15}$ and so $p=x^2+15y^2$,}
\\0&\t{if $p\e 11,14\mod{15}.$}\endcases\endaligned$$
Observe that
$$\f{-15+12\sqrt 5}{105+48\sqrt 5}=(\sqrt 5-2)^2\qtq{and}
\f{42-28\sqrt 5}{-784-350\sqrt 5}=(\sqrt 5-2)^3.$$ Using Lemma 4.1
we see that
$$\align&\sum_{n=0}^{p-1}\Ls{n^3+(-15+12\sqrt 5)n+42-28\sqrt 5}p
\\&=\Ls{\sqrt 5-2}p\sum_{n=0}^{p-1}\Ls{n^3+(105+48\sqrt
5)n-784-350\sqrt 5}p.\endalign$$ Note that $(\sqrt 5-2)\f{1+\sqrt
5}2=\sls{1-\sqrt 5}2^2$. We then have $\sls{\sqrt
5-2}p=\sls{(1+\sqrt 5)/2}p$. Now putting all the above together we
obtain the result.
 \pro{Theorem 4.6} Let $p$ be a prime such that
$p\e 1,4\mod 5$. Then
$$P_{[\f p3]}(\sqrt 5)\e \cases
2x\sls x3\mod p&\t{if $p=x^2+15y^2\e 1,4\mod {15}$,} \\0\mod p&\t{if
$p\e 11,14\mod{15}$}\endcases$$ and
$$\sum_{k=0}^{p-1}\f{\b{2k}k^2\b{3k}k}{(-27)^k}\e
\cases 4x^2\mod p&\t{if $p=x^2+15y^2\e 1,4\mod {15}$,}
\\0\mod {p^2}&\t{if $p\e 11,14\mod {15}$.}
\endcases$$
\endpro
Proof. From Theorem 3.1 and Lemma 4.7 we deduce that
$$\aligned P_{[\f p3]}(\sqrt 5)&\e
-\Ls p3\sum_{n=0}^{p-1}\Ls{n^3+(-15+12\sqrt 5)n+42-28\sqrt 5}p
\\&\e \cases -2x\sls{2x}3=2x\sls x3\mod p&\t{if $p=x^2+15y^2\e
1,4\mod{15}$,}
\\0\mod p&\t{if $p\e 11,14\mod{15}$.}\endcases
\endaligned$$
Taking $m=-27$ and $t=\sqrt 5$ in Theorem 4.1 and then applying the
above we deduce the remaining result. \pro{Conjecture 4.2} Let $p$
be an odd prime such that $p\e 2,7,8,13\mod{15}$. Then
$$P_{[\f p3]}(\sqrt{5})\e \cases 0\mod p&\t{if $p\e
7,13\mod{15}$,}
\\2x\sls x3\sqrt 5\mod p&\t{if $p=5x^2+3y^2\e 2,8\mod {15}$.}
\endcases$$
\endpro
\pro{Lemma 4.8} Let $p$ be a prime such that $p\e \pm 1\mod
 5$. Then
$$\align&\sum_{n=0}^{p-1}\Ls{n^3+(-300+108\sqrt 5)n-2520+1042\sqrt 5}p
\\&=\cases \sls{2\sqrt 5}p\sls x3x&\t{if $p\e 1,4\mod{15}$ and so
$4p=x^2+75y^2$,}\\0&\t{if $p\e 11,14\mod{15}$.}
\endcases\endalign$$
\endpro
Proof. From [I, p.134] we know that the elliptic curve defined by
the equation $y^2=x^3+(-2160+408\sqrt 5)x+42130-10472\sqrt 5$ has
complex multiplication by the order of discriminant $-75$. By (4.4)
and [I, Theorem 3.1] we have
$$\aligned&\sum_{n=0}^{p-1}\Ls{n^3+(-2160+408\sqrt 5)n
+42130-10472\sqrt 5}p
\\&=\cases \sls{-25-13\sqrt 5}p\sls x3x&\t{if $p\e 1,4\mod {15}$
 and so $4p=x^2+75y^2$,}
\\0&\t{if $p\e 11,14\mod{15}.$}\endcases\endaligned$$
Observe that
$$\f{-2160+408\sqrt 5}{-300+108\sqrt 5}=\Big(-\f{7+\sqrt 5}2\Big)^2
\qtq{and} \f{42130-10472\sqrt 5}{-2520+1042\sqrt 5}=\Big(-\f{7+\sqrt
5}2\Big)^3.$$ Using  Lemma 4.1 we see that
$$\align&\sum_{n=0}^{p-1}\Ls{n^3+(-2160+408\sqrt 5)n
+42130-10472\sqrt 5}p
\\&=\Ls{-(7+\sqrt 5)/2}p\sum_{n=0}^{p-1}
\Ls{n^3+(-300+108\sqrt 5)n-2520+1042\sqrt 5}p.\endalign$$ Since
$$2(-7-\sqrt 5)(-25-13\sqrt 5)=120+58\sqrt 5=2\sqrt 5(29+12\sqrt 5)=2\sqrt 5(3+2
\sqrt 5)^2,$$ from the above we deduce the result.

\pro{Theorem 4.7} Let $p$ be a prime such that $p\e 1,4\mod 5$. Then
$$P_{[\f p3]}\Big(\f 9{20}\sqrt 5\Big)\e \cases
-\sls x3x\mod p&\t{if $p\e 1,4\mod {15}$ and so $4p=x^2+75y^2$,}
 \\0\mod p&\t{if
$p\e 11,14\mod{15}$}\endcases$$ and
$$\sum_{k=0}^{p-1}\f{\b{2k}k^2\b{3k}k}{(-8640)^k}\e
\cases x^2\mod p&\t{if $p\e 1,4\mod {15}$ and so $4p=x^2+75y^2$,}
\\0\mod {p^2}&\t{if $p\e 11,14\mod {15}$.}
\endcases$$
\endpro
Proof. From Theorem 3.1 and Lemma 4.8 we deduce that
$$\aligned P_{[\f p3]}\Big(\f 9{20}\sqrt 5\Big)&\e
-\Ls p3\sum_{n=0}^{p-1}\Ls{n^3+\f{27-15\sqrt 5}{\sqrt 5}n
+\f{521-252\sqrt 5}{20}}p
\\&=-\Ls p3\sum_{n=0}^{p-1}\Ls{\sls {n}{2\sqrt 5}^3+\f{27-15\sqrt 5}
{\sqrt 5}\cdot \f {n}{2\sqrt 5} +\f{521-252\sqrt 5}{20}}p
\\&=-\Ls p3\Ls{2\sqrt 5}p\sum_{n=0}^{p-1}
\Ls{n^3+(-300+108\sqrt 5)n-2520+1042\sqrt 5}p
\\&\e \cases -\sls{x}3x\mod p&\t{if $p\e
1,4\mod{15}$ and so $4p=x^2+75y^2$,}
\\0\mod p&\t{if $p\e 11,14\mod{15}$.}\endcases
\endaligned$$
Taking $m=-8640$ and $t=\f 9{4\sqrt 5}$ in Theorem 4.1 and then
applying the above we deduce the remaining result.
\par\q
\newline{\bf Remark 4.4} In [S2] the author conjectured that for any
prime $p>5$,
$$\sum_{k=0}^{p-1}\f{\binom{2k}k^2\binom{3k}k}{(-8640)^k}\e
\cases 4x^2-2p\mod{p^2}&\t{if $3\mid p-1$, $p=x^2+3y^2$ and $5\mid
xy$,}
\\p-2x^2\pm 6xy\mod{p^2}&\t{if $3\mid p-1$, $p=x^2+3y^2$, $5\nmid xy$}
\\&\qq\t{and $x\e \pm y,\pm 2y\mod 5$,}
\\0\mod{p^2}&\t{if $3\mid p-2$.}\endcases$$
This is equivalent to
$$\sum_{k=0}^{p-1}\f{\binom{2k}k^2\binom{3k}k}{(-8640)^k}\e
\cases x^2-2p\mod{p^2}&\t{if $p\e 1,4\mod{15}$ and so
$4p=x^2+75y^2$,}
\\2p-3x^2\mod{p^2}&\t{if $p\e 7,13\mod{15}$ and so
$4p=3x^2+25y^2$,}
\\0\mod{p^2}&\t{if $p\e 2\mod 3$.}\endcases$$

\pro{Conjecture 4.3} For any prime $p>5$ we have
$$\sum_{k=0}^{p-1}\f{9k+1}{(-8640)^k}\b{2k}k^2\b{3k}k\e p\Ls
p{15}\mod{p^3}.$$
\endpro
\pro{Conjecture 4.4} Let $p$ be a prime such that $p\e
7,13,17,23\mod {30}$. Then
$$P_{[\f p3]}\Big(\f 9{20}\sqrt 5\Big)\e \cases
5y\sls y3\mod p&\t{if $p\e 7,13\mod {30}$ and so $4p=3x^2+25y^2$,}
 \\0\mod p&\t{if
$p\e 17,23\mod{30}$.}\endcases$$ \endpro

\par Let $b\in\{17,41,89\}$ and $f(b)=-12^3,-48^3,-300^3$ according
as $b=17,41,89$. In [Su1, Conjectures A20, A22 and A23], Z.W. Sun
conjectured that for any odd prime $p\not=3,b,$
$$\sum_{k=0}^{p-1}\f{\b{2k}k^2\b{3k}k}{f(b)^k}\e
\cases x^2-2p\mod {p^2}&\t{if $\sls p3=\sls pb=1$ and so
$4p=x^2+3by^2$,}
\\2p-3x^2\mod{p^2}&\t{if $\sls p3=\sls pb=-1$ and so
$4p=3x^2+by^2$,}\\0\mod{p^2}&\t{if $\sls p3=-\sls pb$.}
\endcases\tag 4.5$$
\par Now we partially solve (4.5).
\pro{Theorem 4.8} Let $p$ be an odd prime such that $\sls{17}p=1$.
Then
$$P_{[\f p3]}\Ls{\sqrt{17}}4\e \cases  -\sls x3x\mod p&\t{if $p\e 1\mod 3$ and so
$4p=x^2+51y^2$,}
\\0\mod p&\t{if $p\e 2\mod 3$.}\endcases$$ and
$$\sum_{k=0}^{p-1}\f{\b{2k}k^2\b{3k}k}{(-12)^{3k}}
\e\cases x^2\mod p&\t{if $p\e 1\mod 3$ and so $4p=x^2+51y^2$,}
\\0\mod{p^2}&\t{if $p\e 2\mod 3$.}
\endcases$$\endpro
Proof. From [I, p.134] we know that the elliptic curve defined by
the equation $y^2=x^3-(60+12\sqrt {17})x-210-56\sqrt {17}$ has
complex multiplication by the order of discriminant $-51$. Thus, by
(4.4) and [I, Theorem 3.1] we have
$$\aligned&\sum_{n=0}^{p-1}\Ls{n^3-(60+12\sqrt {17})n-210-56\sqrt {17}}p
\\&=\cases \sls {-2}p\sls x3x&\t{if $p\e 1\mod 3$ and so
$4p=x^2+51y^2$,}
\\0&\t{if $p\e 2\mod 3$.}
\endcases\endaligned$$
It then follows from (1.3) and Theorem 3.1 that
$$\aligned P_{[\f p3]}\Ls{\sqrt{17}}4&
=\Ls p3P_{[\f p3]}\Big(-\f{\sqrt{17}}4\Big) \e
-\sum_{n=0}^{p-1}\Ls{n^3+3(-\sqrt{17}-5)n+\f{17}4+22+7\sqrt{17}}p
\\&=-\sum_{n=0}^{p-1}\Ls{(-\f n2)^3-3(5+\sqrt{17})(-\f n2)+\f{105+28\sqrt{17}}
4}p
\\&=-\Ls{-2}p\sum_{n=0}^{p-1}\Ls{n^3-(60+12\sqrt {17})n-210-56\sqrt {17}}p
\\&\e \cases  -\sls x3x\mod p&\t{if $p\e 1\mod 3$ and so
$4p=x^2+51y^2$,}
\\0\mod p&\t{if $p\e 2\mod 3$.}
\endcases\endaligned$$
Taking $m=-12^3$ and $t=\f{\sqrt {17}}4$ in Theorem 4.1 and then
applying the above we deduce the remaining result.

\pro{Conjecture 4.5} Let $p$ be an odd prime such that $\sls
p{17}=-1$. Then
$$P_{[\f p3]}\Ls{\sqrt{17}}4\e \cases 0\mod p&\t{if $p\e 1\mod 3$,}
\\-\sls y3y\sqrt{17}\mod p&\t{if $p\e 2\mod 3$ and so $4p=3x^2+17y^2$.}
\endcases$$
\endpro

\pro{Theorem 4.9} Let $p$ be an odd prime such that $\sls{41}p=1$.
Then
$$P_{[\f p3]}\Ls{5\sqrt{41}}{32}\e \cases  -\sls x3x\mod p&\t{if $p\e 1\mod 3$ and so
$4p=x^2+123y^2$,}
\\0\mod p&\t{if $p\e 2\mod 3$.}\endcases$$ and
$$\sum_{k=0}^{p-1}\f{\b{2k}k^2\b{3k}k}{(-48)^{3k}}
\e\cases x^2\mod p&\t{if $p\e 1\mod 3$ and so $4p=x^2+123y^2$,}
\\0\mod{p^2}&\t{if $p\e 2\mod 3$.}
\endcases$$\endpro
Proof. From [I,p.134] we know that the elliptic curve defined by the
equation $y^2=x^3+(-960+120\sqrt {41})x-13314+2240\sqrt {41}$ has
complex multiplication by the order of discriminant $-123$. Thus, by
[I, Theorem 3.1] we have
$$\aligned&\sum_{n=0}^{p-1}\Ls{n^3-(960-120\sqrt {41})n-13314+2240
\sqrt {41}}p
\\&=\cases \sls {-2}p\sls x3x&\t{if $p\e 1\mod 3$ and so
$4p=x^2+123y^2$,}
\\0&\t{if $p\e 2\mod 3$.}
\endcases\endaligned$$
It then follows from (1.3) and Theorem 3.1 that
$$\aligned P_{[\f p3]}\Ls{5\sqrt{41}}{32}&
 \e
-\Ls p3\sum_{n=0}^{p-1}\Ls{n^3+3(\f
58\sqrt{41}-5)n+\sls{5\sqrt{41}}{16}^2-7\cdot\f{5\sqrt{41}}8+22 }p
\\&=-\Ls p3\sum_{n=0}^{p-1}\Ls{(-\f n8)^3+\f {15}8(\sqrt{41}-8)(-\f n8)
+\f{6657-1120\sqrt{41}}{256}}p
\\&=-\Ls p3\Ls{-2}p\sum_{n=0}^{p-1}\Ls{n^3+(960-120\sqrt{41})n
-13314+2240\sqrt{41}}p
\\&\e \cases  -\sls x3x\mod p&\t{if $p\e 1\mod 3$ and so
$4p=x^2+51y^2$,}
\\0\mod p&\t{if $p\e 2\mod 3$.}
\endcases\endaligned$$
Taking $m=-48^3$ and $t=\f{5\sqrt {41}}{32}$ in Theorem 4.1 and then
applying the above we deduce the remaining result.
 \pro{Conjecture
4.6} Let $p$ be an odd prime such that $\sls p{41}=-1$. Then
$$P_{[\f p3]}\Ls{5\sqrt{41}}{32}\e
\cases 0\mod p&\t{if $p\e 1\mod 3$,}\\ -\sls y3y\sqrt{41}\mod p
&\t{if $p\e 2\mod 3$ and so $4p=3x^2+41y^2$.}\endcases$$
\endpro

\pro{Theorem 4.10} Let $p>5$ be a prime such that $\sls{89}p=1$.
Then
$$P_{[\f p3]}\Ls{53\sqrt{89}}{500}\e
\cases  -\sls x3x\mod p&\t{if $p\e 1\mod 3$ and so $4p=x^2+267y^2$,}
\\0\mod p&\t{if $p\e 2\mod 3$.}\endcases$$ and
$$\sum_{k=0}^{p-1}\f{\b{2k}k^2\b{3k}k}{(-300)^{3k}}
\e\cases x^2\mod p&\t{if $p\e 1\mod 3$ and so $4p=x^2+267y^2$,}
\\0\mod{p^2}&\t{if $p\e 2\mod 3$.}
\endcases$$\endpro
Proof. From [I, p.135] we know that the elliptic curve defined by
the equation $y^2=x^3+(-37500+3180\sqrt {89})x+3250002-371000\sqrt
{89}$ has complex multiplication by the order of discriminant
$-267$. Thus, by [I, Theorem 3.1] we have
$$\aligned&\sum_{n=0}^{p-1}\Ls{n^3+(-37500+3180\sqrt {89})n
+3250002-371000\sqrt {89}}p
\\&=\cases \sls {2}p\sls x3x&\t{if $p\e 1\mod 3$ and so
$4p=x^2+267y^2$,}
\\0&\t{if $p\e 2\mod 3$.}
\endcases\endaligned$$
It then follows from (1.3) and Theorem 3.1 that
$$\aligned P_{[\f p3]}\Ls{53\sqrt{89}}{500}&
 \e
-\Ls p3\sum_{n=0}^{p-1}\Ls{n^3+3(4\cdot\f{53}{500}\sqrt{89}-5)n+
4(\f{53}{500}\sqrt{89})^2-28\cdot\f{53\sqrt{89}}{500}+22 }p
\\&=-\Ls p3\sum_{n=0}^{p-1}\Ls{(\f n{50})^3
+\f{-1875+159\sqrt{89}}{125}\cdot \f
n{50}+\f{1625001-185500\sqrt{89}}{250^2}}p
\\&=-\Ls p3\Ls{50}p\sum_{n=0}^{p-1}\Ls{n^3+(-37500+3180\sqrt {89})n
+3250002-371000\sqrt {89}}p
\\&\e \cases  -\sls x3x\mod p&\t{if $p\e 1\mod 3$ and so
$4p=x^2+51y^2$,}
\\0\mod p&\t{if $p\e 2\mod 3$.}
\endcases\endaligned$$
Taking $m=-300^3$ and $t=\f{53\sqrt {89}}{500}$ in Theorem 4.1 and
then applying the above we deduce the remaining result.
\pro{Conjecture 4.7} Let $p$ be an odd prime such that $\sls
p{89}=-1$. Then
$$P_{[\f p3]}\Ls{53\sqrt{89}}{500}\e \cases 0\mod p&\t{if $p\e 1\mod 3$,}
\\-\sls y3y\sqrt{89}\mod p
&\t{if $p\e 2\mod 3$ and so $4p=3x^2+89y^2$.}\endcases$$
\endpro

\par In the end we pose the following conjectures.
\pro{Conjecture 4.8} Let $p$ be a prime with $p\e \pm 1\mod 5$. Then
$$\aligned&\sum_{n=0}^{p-1}\Ls{n^3+3(-125+44\sqrt 5)n+154(21-10\sqrt 5}p
\\&=\cases 2x\sls {2x}3&\t{if $p\e 1,4\mod {15}$ and so $p=x^2+15y^2$,}
\\0&\t{if $p\e 11,14\mod{15}.$}\endcases\endaligned$$
\endpro
If Conjecture 4.8 is true, we may discuss the following conjecture
in [S2]:
$$\sum_{k=0}^{p-1}\f{\binom{2k}k^2\binom{3k}k}{15^{3k}}\e\cases 4x^2-2p\mod {p^2}
&\t{if $p\e 1,4\mod {15}$ and so $p=x^2+15y^2$,}
\\2p-12x^2\mod {p^2}
&\t{if $p\e 2,8\mod {15}$ and so $p=3x^2+5y^2$,}
\\ 0\mod{p^2}&\t{if $p\e 7,11,13,14\mod {15}$.}
\endcases$$
\pro{Conjecture 4.9} For any prime $p>5$ we have
$$\sum_{k=0}^{p-1}\f{33k+4}{15^{3k}}\b{2k}k^2\b{3k}k\e
4p\Ls p3\mod{p^3}.$$
\endpro
\pro{Conjecture 4.10} Let $p>5$ be a prime. Then
$$P_{[\f p3]}\Ls{11\sqrt 5}{25}
\e\cases 2x\sls x3\mod p&\t{if $p\e 1,19\mod{30}$ and so
$p=x^2+15y^2$,}
\\-2y\sls y3\sqrt 5\mod p&\t{if $p\e 17,23\mod {30}$ and so
$p=3x^2+5y^2$,}
\\0\mod p&\t{if $p\e 7,11,13,29\mod{30}$.}
\endcases$$
\endpro

\par In [S2], the author also conjectured that for any prime $p>3$,
$$\sum_{k=0}^{p-1}\f{\binom{2k}k^2\binom{3k}k}{1458^k}\e\cases 4x^2-2p\mod {p^2}
&\t{if $p\e 1\mod 3$ and so $p=x^2+3y^2$,}
\\ 0\mod{p^2}&\t{if $p\e 2\mod 3$.}
\endcases$$

\pro{Conjecture 4.11} For any prime $p>3$ we have
$$\sum_{k=0}^{p-1}\f{15k+2}{1458^k}\b{2k}k^2\b{3k}k\e (-1)^{\f{p-1}2}2p
\mod{p^3}.$$
\endpro

\pro{Conjecture 4.12} Let $p>3$ be a prime. Then
$$P_{[\f p3]}\Ls{5\sqrt 3}9
\e\cases (-1)^{\f{p-1}2}2x\sls x3\mod p&\t{if $p\e 1\mod 3$ and so
$p=x^2+3y^2$,}
\\0\mod p&\t{if $p\e 2\mod 3$.}
\endcases$$
\endpro

\enddocument
\bye